\documentclass[12pt, twoside, a4paper, fleqn]{article}

\usepackage{latexsym}
\usepackage{amscd}
\usepackage{theorem}
\usepackage{pifont}
\usepackage{mathbbol}
\usepackage{amsfonts}
\usepackage{xspace}
\usepackage{amssymb}
\usepackage{fancyhdr}
\usepackage{mathrsfs}
\usepackage{amsmath}%
\usepackage{graphics}%
\usepackage{graphicx}%
\usepackage{euscript}%
\usepackage{psfrag}%
\usepackage{empheq}%
\usepackage{upgreek}
\usepackage{eufrak}

\DeclareMathAlphabet{\mathpzc}{OT1}{pzc}{m}{it}

{\theorembodyfont{\slshape} \newtheorem{theorem}{Theorem}[section]}
{\theorembodyfont{\slshape} \newtheorem{lemma}[theorem]{Lemma}}
{\theorembodyfont{\slshape} }
{\theorembodyfont{\slshape} }
\newcommand{\norm}[1]{\left\lVert#1\right\rVert}

\newenvironment{proof}{\noindent\textbf{Proof:\ }}{$\hfill{\bullet}$}

\numberwithin{equation}{section}

\title{\textsc{Distribution of Typical Orbits for a Skew-Product Map Generated by Random Dynamics of Finitely many Rational Maps}}                

\author{Shrihari Sridharan \\ {\tt shrihari@iisertvm.ac.in}  \bigskip \\ Sharvari Neetin Tikekar \\ {\tt sharvai.tikekar14@iisertvm.ac.in} \bigskip \\ Atma Ram Tiwari \\ {\tt artiwari15@iisertvm.ac.in} \bigskip 
 \bigskip \\ 
 Indian Institute of Science Education and Research \\ Thiruvananthapuram(IISER-TVM), India.}     

\date{\today}

\begin{document}

\maketitle
\thispagestyle{empty}

\begin{abstract} 
\noindent 
In this paper, we consider the dynamics of a skew-product map defined on the Cartesian product of the symbolic one-sided shift space on $N$ symbols and the complex sphere where we allow $N$ rational maps, $R_{1}, R_{2}, \cdots, R_{N}$, each with degree $d_{i};\ 1 \le i \le N$ and with at least one $R_{i}$ in the collection whose degree is at least $2$. We obtain results regarding the distribution of pre-images of points and the periodic points in a subset of the product space (where the skew-product map does not behave normally). We further explore the ergodicity of the Sumi-Urbanskii (equilibrium) measure associated to some real-valued H\"{o}lder continuous function defined on the Julia set of the skew-product map and obtain estimates on the mean deviation of the behaviour of typical orbits, violating such ergodic necessities. 
\end{abstract}
\bigskip \bigskip

\begin{tabular}{l c l}
\textbf{Keywords} & : & Equidistribution (Lyubich's) measure for \\ & & skew-product maps, \\ & & Sumi-Urbanskii equilibrium measure, \\ & & Large deviation results. \\
\\
\textbf{AMS Subject Classifications} & : & 30D05, 32H50, 37F10.
\end{tabular}
\bigskip 
\bigskip 

\section{Introduction}

\noindent 
Investigations into the theory of the dynamics of semigroups of rational functions were initiated by Hinkkanen and Martin in \cite{hm:96:1, hm:96:2, hm:96:3}, and a lot is known, by now about such dynamical systems. Around the same time, Zhou and Ren \cite{rz:92} and Gong and Ren \cite{gr:96} introduced a different approach to study random iteration of rational maps. The study of the skew-product map, initiated by Sumi in \cite{hs:00, hs:01, hs:06} is a very useful tool to study the dynamics of finitely many rational maps acting on the Riemann sphere, in tandem. In fact, skew-product maps are natural objects that one finds in random dynamics with finitely many rational maps, equipped with the necessary book-keeping tools to maintain track of orbits. The first of the parameters helps us in identifying the rational map to be put to use out of the finitely many rational maps that we employ to act on the Riemann sphere. A natural dichotomy of the Riemann sphere into Fatou and Julia sets based on the normality of a growing random sequence of maps (in the sense of Montel) was achieved by Sumi. This situation generalises the setting where we consider a single rational map with degree at least $2$ and dichotomise the Riemann sphere based on the normality of the sequence of iterates (also, in the sense of Montel). Here and throughout this paper, we mean by a rational map that it can be written as the quotient of any two relatively prime polynomials, \textit{i.e.}, in its irreducible form. Hence, it is natural to investigate the veracity of the various properties prevalent in the latter setting of the dynamics under a single rational map to the former setting of random dynamics under finitely many rational maps. 
\medskip 

\noindent 
One of the important properties that the dynamics of a single polynomial map enjoys is the accumulation of pre-images of any generic point in the Riemann sphere over a compact subset called the Julia set of the Riemann sphere, where the family of iterates of the polynomial map is not normal. In fact, Brolin proves in \cite{hb:65} that the normalised sum of point masses supported on the set of pre-images for any generic point in the Riemann sphere under any polynomial map defined there is equidistributed. Lyubich considered generic points not from the whole Riemann sphere, but from the Julia set however for a rational map and obtained more taut results concerning such distribution of points in the Julia set, in \cite{myl:82, myl:86}. The results due to Lyubich come in two varieties, one exploiting the density of the pre-images of any arbitrary point in the rational map-invariant Julia set, while the other exploiting the density of periodic points of the rational map in the Julia set. 
\medskip 

\noindent 
For such a rational map $R$ restricted on its Julia set $J_{R}$, we have a valid theory of thermodynamic formalism due to authors like Bowen, Ruelle, Coelho and Parry, Parry and Pollicott, Zinsmeister etc. \cite{rb:72, dr:78, cp:90, pp:90, mz:00}. A particularly desirable feature of this analysis is the definition of the pressure of any real-valued continuous function, $f$ defined on the Julia set $J_{R}$. Suppose the real-valued continuous function $f$ further satisfies a H\"{o}lder condition for some H\"{o}lder exponent $0 < \alpha < 1$, then Denker and Urbanski prove the unique existence of an equilibrium $R$-invariant non-atomic probability measure, say $\mu_{f}$ supported on $J_{R}$ in \cite{du:91, dpu:96}. 
\medskip 

\noindent 
Based on the equidistribution of pre-images of a generic point in $J_{R}$ due to Lyubich and considering a probability orbital measure, \textit{i.e.,} a probability measure defined on the orbit of an arbitrary point chosen from a dense subset of $J_{R}$, Pollicott and Sharp, \cite{ps:96} obtained a large derivation result by verifying the rate at which the orbital measure does not have a membership in an arbitrary weak*-open neighbourhood of the equilibrium measure $\mu_{f}$. This result describes the ergodicity of typical orbits under the setting. In fact, the authors estimate and obtain the exponential speed at which any generic orbit decays, if burdened with a condition of non-ergodicity, as described above. An analogous result in the setting of periodic points was obtained by Pollicott and Sridharan in \cite{ps:07}. 
\medskip 

\noindent 
Our aim in this paper is to generalise the results explained above due to Lyubich \cite{myl:86}, Pollicott and Sharp \cite{ps:96} and Pollicott and Sridharan \cite{ps:07} to the setting of dynamics under finitely many rational maps acting on the Riemann sphere. The reader will be aware of a result due to Boyd in \cite{db:99} that may seem to address one aspect of the Lyubich's result. However, we generalise the scope of Boyd's result as well, by relaxing the condition that every rational map in the collection must be of degree at least $2$. We shall contend with the condition that there exists at least one rational map in the collection whose degree is at least $2$. 
\medskip 

\noindent 
The paper is organised as follows. In the next section \eqref{bssr}, we build the necessary mathematical objects and explain the basic setting in which the main results of this paper are true. We shall give the essential definitions and write the precise statements of the theorems as well. Once the main theorems are stated, we will prove them one by one in the subsequent sections. In section \eqref{techcond}, we define two different operators on the Banach space of real-valued continuous functions defined on the Julia set and study certain properties of those operators. Then we will state and prove some important technical lemmas in section \eqref{ttl}. These technical lemmas will come in handy to write the proof of the theorem \eqref{edpr}, in section \eqref{pfpt2}. In sections \eqref{pfper} and \eqref{pfld}, we complete the proofs of the rest of the theorems \eqref{lyuper}, \eqref{lard-preimage} and \eqref{periodicpts}. 

\section{Basic setting and the statements of results} 
\label{bssr}

\noindent 
In this section, we shall narrate the elementary setting to that extent that we can state our main results of this paper. 
\medskip

\noindent 
Let $\overline{\mathbb{C}} $ denote the Riemann sphere, meaning the complex plane along with the point at $\infty$. Let $\left\{ R_{1}, R_{2}, \cdots, R_{N} \right\}$ be a finite collection of rational maps defined on $\overline{\mathbb{C}}$, each with degree $d_{i};\ 1 \le i \le N$ respectively and at least one $R_{i}$ in the collection whose degree is at least $2$. We now introduce the symbolic space of one-sided infinite sequences on $N$ symbols, $\Sigma_{N}^{+} := \left\{ 1, 2, \cdots, N \right\}^{\mathbb{Z}_{+}}$ that helps us in keeping track of the map that we choose from the collection of rational maps, as described below in the skew-product map. 
\medskip

\noindent 
Define the \emph{skew-product map} $T : X := \Sigma_{N}^{+} \times \overline{\mathbb{C}} \longrightarrow X$ as 
\begin{equation} 
T\left(\left(w, z\right)\right)\ \ :=\ \ \left( \sigma w, R_{v_{1}} z\right), 
\end{equation} 
where $w = (v_{1}\, v_{2}\, \cdots )$ with $v_{i} \in \left\{ 1, 2, \cdots, N \right\}$ and $\sigma : \Sigma_{N}^{+} \longrightarrow \Sigma_{N}^{+}$ is the usual shift map that shifts the sequences one place to the left, $(\sigma w)_{n} = v_{n + 1}$ for $n \ge 1$. Throughout this paper, we will reserve $v_{i}$ to denote letters from $\left\{ 1, 2, \cdots, N \right\}$ and $w_{i}$ to denote the infinite-lettered words from $\Sigma_{N}^{+}$. Observe that the action of $T$ on the second parameter can also be described as a semi-group action on $\overline{\mathbb{C}}$ generated by the collection of rational maps, $G := \langle R_{1}, R_{2}, \cdots, R_{N} \rangle$ with composition of maps acting as the group operation.  We denote by $\Pi_{\Sigma_{N}^{+}} : X \longrightarrow \Sigma_{N}^{+}$ and $\Pi_{\overline{\mathbb{C}}} : X \longrightarrow \overline{\mathbb{C}}$, the projection maps of $X$ onto $\Sigma_{N}^{+}$ and $\overline{\mathbb{C}}$ respectively. Note that under the canonical identification, $\left( \Pi_{\Sigma_{N}^{+}} \right)^{-1} \left( \{ w \} \right) \cong \overline{\mathbb{C}}  $, \textit{i.e.,} each fiber $ \Pi_{\Sigma_{N}^{+}}^{-1} \left( \{ w \} \right)$ is isomorphic to the Riemann sphere, $\overline{\mathbb{C}}$. 
\medskip 

\noindent 
We shall consider the spherical metric $d_{\overline{\mathbb{C}}} (\cdot, \cdot)$ on $\overline{\mathbb{C}}$, while the metric $d_{\Sigma_{N}^{+}}$ on the symbolic space is defined by $d_{\Sigma_{N}^{+}} (w_{1}, w_{2}) := \frac{1}{2^{n(w_{1},\, w_{2})}}$ where $n(w_{1}, w_{2})$ is the first place where the infinitely long words $w_{1}$ and $w_{2}$ disagree. We set $n(w, w) := \infty$, in order that $d_{\Sigma_{N}^{+}} (w, w) = 0$. We now define the product metric on $X$ as 
\[ d_{X} \left( \left( w_{1}, z_{1} \right), \left( w_{2}, z_{2} \right) \right)\ \ :=\ \ \max \left\{d_{\Sigma_{N}^{+}} (w_{1}, w_{2}), d_{\overline{\mathbb{C}}} (z_{1}, z_{2})\right\}. \] 

\noindent
For $w  = (v_{1}\, v_{2}\, \cdots ) \in \Sigma_{N}^{+}$, let $F_{w}$ be the collection of points $z \in \overline{\mathbb{C}}$ such that the family 
\[ \left\{ R_{v_{n}} \circ R_{v_{n - 1}} \circ \cdots \circ R_{v_{1}} \right\}_{n \ge 1} \] 
is normal (in the sense of Montel) in a neighbourhood of $z$ and $J_{w} := \overline{\mathbb{C}} \setminus F_{w}$. We then define the \emph{Julia set} $\widetilde{ \mathcal{J} }$ and the \emph{Fatou set} $\widetilde{ \mathcal{F} }$ respectively of the skew-product map to be  
\[ \widetilde{ \mathcal{J} }\ :=\ \overline{\bigcup_{w\, \in\, \Sigma_{N}^{+}} \left\{w\right\} \times J_{w}};\ \ \ \ \widetilde{ \mathcal{F} }\ :=\ X \setminus \widetilde{ \mathcal{J} }. \]
The closure in the definition of $\widetilde{ \mathcal{J} }$ is taken with respect to the product metric as described above. Then, it is obvious that $\widetilde{ \mathcal{J} }$ is a compact metric space. In fact, the Julia set of the skew-product map, $\widetilde{ \mathcal{J} }$ satisfies the usual properties of the Julia set of a rational map. For more properties of the Fatou set and the Julia set of the skew-product map, interested readers are referred to \cite{hs:00}. 
\medskip 

\noindent 
For any set $E \subseteq X$, we use the notations $E_{\Sigma_{N}^{+}} = \Pi_{\Sigma_{N}^{+}} (E)$ and $E_{\overline{\mathbb{C}}} = \Pi_{\overline{\mathbb{C}}} (E)$ to be the projections of $E$ onto $\Sigma_{N}^{+}$ and $\overline{\mathbb{C}}$ respectively, throughout this paper. We say a set $E \subset X$ is \emph{simply connected} if $E_{\overline{\mathbb{C}}}$ is simply connected. Further, the \textit{boundry} of the set $E$ is defined as $\partial E := E_{\Sigma_{N}^{+}} \times \partial E_{\overline{\mathbb{C}}} $. Note that these definitions makes total sense, since the space $\Sigma_{N}^{+}$ is totally disconnected and we use it only for book-keeping purposes. 
\medskip 

\noindent
Consider the sequence of probability mass distributions $\mu_{n}^{((w_{0},\, z_{0}))}$ supported on the $n^{\rm th}$ order preimages of an arbitrarily chosen  point $(w_{0}, z_{0}) \in X$, \textit{i.e.,}
\begin{equation} 
\label{mass_dis_preim_eq}
\mu_{n}^{((w_{0},\, z_{0}))}\ \ :=\ \ \frac{1}{\left( d_{1} + d_{2} + \cdots + d_{N} \right)^{n}} \sum_{T^{n} ((w,\, z)) = (w_{0},\, z_{0})} \delta_{(w,\, z)},
\end{equation}
where $\delta_{(w,\, z)}$ is the Dirac delta mass concentrated at the point $(w, z)$. Observe that the constant factor in the definition of the measure is used precisely for the purpose of rendering the quantity to be a probability for every $n$. We denote the support of any measure $\mu$ on $X$ by ${\rm supp} \left( \mu \right)$.
\medskip

\noindent 
\begin{theorem} 
\label{edpr} 
For a generic point $(w_{0}, z_{0}) \in \widetilde{ \mathcal{J} }$, the sequence $\left\{ \mu_{n}^{((w_{0},\, z_{0}))} \right\}_{n \ge 1}$ converges to some probability measure $\mu^{*}$ (in the weak*-topology), independent of $(w_{0}, z_{0})$. Also, ${\rm supp} \left( \mu^{*} \right) = \widetilde{\mathcal{J}}$.  
\end{theorem}
\medskip 

\noindent 
The above result is a generalisation of the Lyubich's theorem in the case of a rational map of degree at least $2$ and can be compared to Boyd's result, as explained in the introduction. As a concerted reader may know, Lyubich's results on distribution of points in the case of a rational map of degree at least $2$ comes in two flavours; one concerning the equidistribution of pre-images of a generic point in the Julia set and the other concerning the equidistribution of the periodic points of the rational map. In our context, a point $(\omega, \zeta) \in X$ is called a \emph{periodic point} of period $m$ for $T$ if $T^{m} (\omega, \zeta) = (\omega, \zeta)$ \textit{i.e.,} $\sigma^{m} \omega = \omega$ and $\left( R_{\vartheta_{m}} \circ R_{\vartheta_{m - 1}} \circ \cdots \circ R_{\vartheta_{1}} \right) (\zeta) = \zeta$, where $\omega = (\vartheta_{1}\, \vartheta_{2}\, \cdots\, \vartheta_{m}\, \vartheta_{1}\, \cdots\, \vartheta_{m}\, \cdots)$. We shall now state the analogous result in the case of periodic points. 
\medskip 

\noindent 
\begin{theorem} 
\label{lyuper}
Consider the sequence of probability mass distributions $\nu_{m}$ supported on the periodic points of order $m$ in $X$, \textit{i.e.},
\begin{equation} 
\nu_{m} := \frac{1}{(d_{1} + d_{2} + \cdots + d_{N})^{m} + N} \sum_{(\omega,\, \zeta)\, \in\, {\rm Fix}_{m}(T)} \delta_{(\omega,\, \zeta)}, 
\end{equation} 
where ${\rm Fix}_{m} (T)$ is the collection of all fixed points of $T^{m}$. Then, the sequence $\left\{ \nu_{m} \right\}_{m \ge 1}$ converges (in the weak*-topology) to $\mu^{*}$.
\end{theorem}
\medskip 

\noindent 
The measure $\mu^{*}$ that one obtains using any of the above mentioned two procedures for the compact Julia set $\widetilde{\mathcal{J}}$ is a non-atomic, $T$-invariant, ergodic measure that is uniformly distributed on its support, namely $\widetilde{ \mathcal{J} }$, with no point carrying any mass. 
\medskip 

\noindent 
We shall, henceforth focus on writing other necessary definitions and some results from the literature that will, in turn help us in stating certain large deviation results pertaining to the pre-images of a generic point as well as periodic points. Suppose $\mathcal{M}_{T}$ denotes the space of all $T$-invariant probability measures defined on $\widetilde{\mathcal{J}}$. For any continuous function $f: \widetilde{\mathcal{J}} \longrightarrow \mathbb{R}$, the \emph{pressure} of $f$ is defined as, 
\[ \mathfrak{P} (f)\ \ :=\ \ \sup_{m\, \in\, \mathcal{M}_{T}} \left\{ h_{m} (T) + \int_{\widetilde{\mathcal{J}}} f d m \right\}, \] 
where $h_{m} (T)$ is the entropy of the skew-product map $T$ with respect to the measure $m$. Observe that when $f \equiv 0,\ \mathfrak{P} (0)$ evaluates the topological entropy of $T$ restricted on its Julia set, $\widetilde{\mathcal{J}}$, given by 
\[ \mathfrak{P} (0)\ \ =\ \ \sup_{m\, \in\, \mathcal{M}_{T}} h_{m} (T)\ \ =\ \ \log \left( \sum_{i = 1}^{N} d_{i} \right). \] 
In fact, upon satisfying a further pair of technical conditions (that we shall explain in section \eqref{techcond}), Sumi and Urbanski, in \cite{su:09}, proves that there exists a unique equilibrium state denoted by $\mu_{f}$ for every real-valued H\"{o}lder continuous function $f$. 
\medskip 

\noindent 
We conclude this section with two more theorems concerning the deviation from mean ergodic behaviour of typical orbits (periodic or otherwise) with respect to the equilibrium measures for real-valued H\"{o}lder continuous functions. We begin with the following definitions of orbital measures. For any generic point $(w_{0}, z_{0}) \in \widetilde{ \mathcal{J} }$ with $(w_{n}, z_{n}) \in T^{-n} \left((w_{0}, z_{0})\right)$, define a sequence of normalised proportion of orbital measure as  
\begin{equation} 
\label{one}
\mu_{(w_{n},\, z_{n})}\ :=\ \frac{1}{n} \left[ \delta_{(w_{n},\, z_{n})} + \delta_{T \left((w_{n},\, z_{n})\right)} + \cdots + \delta_{T^{n - 1} \left((w_{n},\, z_{n})\right)} \right]. 
\end{equation} 
And for a periodic point $(\omega, \zeta) \in {\rm Fix}_{m} (T)$, the normalised proportion of the periodic orbital measure is defined as  
\begin{equation} 
\label{two}
\nu_{(\omega,\, \zeta),\, m}\ :=\ \frac{1}{m} \left[ \delta_{(\omega,\, \zeta)} + \delta_{T \left((\omega,\, \zeta)\right)} + \cdots + \delta_{T^{m - 1} \left((\omega,\, \zeta)\right)} \right]. 
\end{equation} 

\noindent 
\begin{theorem} 
\label{lard-preimage} 
Let $f : \widetilde{\mathcal{J}} \longrightarrow \mathbb{R}$ be a H\"{o}lder continuous function with equilibrium measure $\mu_{f}$ that satisfies $\mathfrak{P}(f) > \sup f + \log N$. Suppose there exists some weak*-open neighbourhood $\mathcal{U}$ around $\mu_{f}$ in the space $\mathcal{M}_{T}$ such that all but finitely many elements of the sequence of the normalised proportion of some orbital measure $\mu_{(w_{n},\, z_{n})} \notin \mathcal{U}$. Then, the sequence of measures $\left\{ \mu_{(w_{n},\, z_{n})} \right\}_{n \ge 1}$ (as defined in \eqref{one}) converges to zero exponentially fast. 
\end{theorem} 
\medskip 

\noindent 
\begin{theorem} 
\label{periodicpts}
Let $f$ and $\mu_{f}$ be as in the theorem \eqref{lard-preimage}. Suppose there exists some weak*-open neighbourhood $\mathcal{V}$ around $\mu_{f}$ in the space $\mathcal{M}_{T}$ such that all but finitely many elements of the sequence of the normalised proportion of the periodic orbital measure $\nu_{(\omega,\, \zeta), m} \notin \mathcal{V}$, where $(\omega, \zeta) \in \widetilde{\mathcal{J}} \cap {\rm Fix}_{m} (T)$. Then, the sequence of measures $\left\{ \nu_{(\omega,\, \zeta),\, m} \right\}_{m \ge 1}$ (as defined in \eqref{two}) converges to zero exponentially fast. 
\end{theorem} 

\section{Two operators on $\mathcal{C} \left( \widetilde{\mathcal{J}}, \mathbb{R} \right)$} 
\label{techcond}

\noindent 
In this section we define two different operators on the Banach space of real-valued continuous functions defined on $\widetilde{\mathcal{J}}$.  
\medskip 

\noindent 
For $w = (v_{1}\, v_{2}\, \cdots ) \in \Sigma_{N}^{+},\ (w, z) \in X$ and $m \in \mathbb{Z}_{+}$, we set 
\[ (T^{m})'((w, z)) := \left( R_{v_{m}} \circ R_{v_{m - 1}} \circ \cdots \circ R_{v_1} \right)' (z). \] 
A point $p = (w, z) \in X$ for $w = \left( v_{1}\, v_{2}\, \cdots \right)$ is called a \emph{critical point} for the skew-product map $T$ if $\Pi_{\overline{\mathbb{C}}} (p)$ is a critical point of $R_{v_{1}}$. We denote the set of all critical points of $T$ by $\mathfrak{C}$. Further, a periodic point $(\omega, \zeta) \in X$ of period $m$ for $T$ is classified as \emph{super-attracting} if $|(T^{m})'((\omega, \zeta))| = 0$. We denote by $(\omega_{sa}, \zeta_{sa})$, a super-attracting fixed point of $T$. 
\medskip

\noindent
Throughout this paper, we shall demand that the Julia set $\widetilde{\mathcal{J}}$ remains free of the super-attracting fixed points of $T$, \textit{i.e.,} 
\begin{enumerate} 
\item[{\bf (A1)}] $(\omega_{sa}, \zeta_{sa})\ \ \notin\ \ \widetilde{\mathcal{J}}$. 
\end{enumerate} 
\medskip

\noindent 
Denote by $\mathcal{C} \left( \widetilde{\mathcal{J}}, \mathbb{R} \right)$, the Banach space of all real-valued continuous functions defined on $\widetilde{\mathcal{J}}$ equipped with the usual supremum norm given by $\norm{f} = \sup_{p\, \in\, \widetilde{\mathcal{J}}} | f(p) |$. Further, for any $f \in \mathcal{C} \left( \widetilde{\mathcal{J}}, \mathbb{R} \right)$, we call by $f^{n}(p)$, the  ergodic sum of $f$ of order $n$, that is,
\begin{eqnarray*} 
f^{n} (p)\ =\ f^{n} ((w, z)) & = & f((w, z)) + f(T (w, z)) + \cdots + f(T^{n - 1} (w, z)) \\ 
& = & f((w, z)) + f((\sigma w, R_{v_{1}} z)) + \cdots + f((\sigma^{n - 1} w, R_{v_{n - 1}} \circ \cdots \circ R_{v_{1}} z)). 
\end{eqnarray*} 

\noindent 
Consider the \emph{Perron-Frobenius operator} $\mathcal{L}_{f}$ for some real-valued H\"{o}lder continuous function $f \in \mathcal{C}^{\alpha} \left( \widetilde{\mathcal{J}}, \mathbb{R} \right)$ defined on $\mathcal{C} \left( \widetilde{\mathcal{J}}, \mathbb{R} \right)$ by 
\[ \left( \mathcal{L}_{f} g \right) ((w, z))\ \ :=\ \ \sum_{T\left((w_{1},\, z_{1})\right)\, =\, (w,\, z)} e^{f((w_{1},\, z_{1}))} g((w_{1}, z_{1})). \] 
Observe that this definition of the Perron-Frobenius operator then entails, 
\[ \left( \mathcal{L}_{f}^{n} g \right) ((w, z))\ \ = \sum_{T^{n}(w_{n},\, z_{n})\, =\, (w,\, z)} e^{f^{n}((w_{n},\, z_{n}))}  g((w_{n}, z_{n})). \] 
Further, for any H\"{o}lder continuous function $f \in \mathcal{C}^{\alpha} \left( \widetilde{\mathcal{J}}, \mathbb{R} \right)$, we define the \emph{pointwise pressure} denoted by $\mathfrak{P}^{p} (f)$ as 
\begin{equation} 
\label{ppr}
\mathfrak{P}^{p} (f)\ \ :=\ \ \sup \left\{ P_{(w,\, z)} (f) : (w, z) \in \widetilde{ \mathcal{J} } \right\},\ \ \text{where}\ \ P_{(w,\, z)} (f)\ :=\ \limsup_{n \to \infty} \frac{1}{n} \log \left( \mathcal{L}_{f}^{n} \mathbf{1} \right) \left((w, z)\right). 
\end{equation}  
Then, a theorem due to Sumi and Urbanski states as follows. 
\medskip 

\noindent 
\begin{theorem}\cite{su:09} 
\label{SumiUrbanski}
Let $f \in \mathcal{C}^{\alpha} \left( \widetilde{\mathcal{J}}, \mathbb{R} \right)$ such that $\mathfrak{P}^{p} (f) > \sup f + \log N$. Then, there exists a unique equilibrium state $\mu_{f} \in \mathcal{M}_{T}$ for the function $f$. Moreover, the supremum is achieved for every point $(w, z) \in \widetilde{\mathcal{J}}$, in the definition of pointwise pressure, \eqref{ppr}. 
\end{theorem}
\medskip 

\noindent 
Thus, this theorem gives us an alternate definition of the pressure function, namely, 
\[ \mathfrak{P} (f)\ \ =\ \ \mathfrak{P}^{p} (f)\ \ =\ \ \lim_{n \to \infty} \frac{1}{n} \log \left( \mathcal{L}_{f}^{n} \mathbf{1} \right) \left((w, z)\right), \] 
for any generic $(w, z) \in \widetilde{\mathcal{J}}$. In fact, owing to the density of pre-images of generic points in $\widetilde{\mathcal{J}}$ and the density of periodic points in $\widetilde{\mathcal{J}}$, one can as well observe that 
\[ \lim_{n \to \infty} \frac{1}{n} \log \sum_{T^{n} \left((w_{n},\, z_{n})\right)\, =\, \left((w,\, z)\right)} e^{f^{n} \left((w_{n},\, z_{n})\right)}\ \ =\ \ \lim_{m \to \infty} \frac{1}{m} \log \sum_{(\omega,\, \zeta)\, \in\, {\rm Fix}_{m} (T)} e^{f^{m} \left((\omega,\, \zeta)\right)}. \] 
We now recall a result due to Pollicott and Sharp, as in \cite{ps:96}. 
\medskip 

\noindent 
\begin{theorem}\cite{ps:96} 
\label{ps:96} 
Whenever $f \in \mathcal{C}^{\alpha} \left( \widetilde{\mathcal{J}}, \mathbb{R} \right)$ and $g \in \mathcal{C} \left( \widetilde{\mathcal{J}}, \mathbb{R} \right)$, we have 
\[ \limsup_{n \to \infty} \frac{1}{n} \log \left( \mathcal{L}_{f + g}^{n} \mathbf{1} \right) \left((w, z)\right)\ \ \le\ \ \mathfrak{P} (f + g). \] 
\end{theorem}
\medskip 

\noindent 
The second operator that we define on $\mathcal{C} \left( \widetilde{\mathcal{J}}, \mathbb{R} \right)$ is a continuous linear operator $\mathcal{G}$ defined by 
\begin{eqnarray} 
\left( \mathcal{G} h \right) \left((w, z)\right) & := & \int_{\widetilde{\mathcal{J}}} h \left((w_{1}, z_{1})\right) d \mu_{1}^{\left((w,\, z)\right)} \left((w_{1}, z_{1})\right) \\ 
& = & \frac{1}{d_{1} + d_{2} + \cdots + d_{N}} \sum_{v_{1} = 1}^{N} \sum_{j = 1}^{d_{v_{1}}} h \left(((v_{1}\, w), z_{(v_{1},\, j)})\right), 
\end{eqnarray} 
where $((v_{1}\, w), z_{(v_{1},\, j)})$ are the solutions of the equation $T \left(((v_{1}\, w), z_{(v_{1},\, j)})\right) = (w, z)$. Extending the above notation for the solution set, we remark that for any point $p \in X$, the concatenation of all possible $n$-lettered words, say $(v_{n}\, v_{n - 1}\, \cdots\, v_{1})$ to $\Pi_{\Sigma_{N}^{+}} (p)$ gives the set of all $n$-th order pre-images on the first parameter of $X$ while the set 
\[ \bigcup_{v_{n} = 1}^{N} \bigcup_{v_{n - 1} = 1}^{N} \cdots \bigcup_{v_{1} = 1}^{N} \bigcup_{j = 1}^{d_{v_{n}} d_{v_{n - 1}} \cdots d_{v_{1}}} \left\{ z_{((v_{n}\, v_{n - 1}\, \cdots\, v_{1}),\, j)} \right\}, \] 
gives the set of all $n$-th order pre-images with appropriate labels on the second parameter of $X$. Observe that the cardinality of the above mentioned set is $(d_{1} + d_{2} + \cdots + d_{N})^{n}$. We also note that as $(w, z)$ moves continuously in $X$, the solution set $\left\{ ((v_{1}\, w), z_{(v_{1},\, j)}) \right\}$ also moves continuously in $X$. Hence, $\mathcal{G} h$ can be treated as a continuous function defined on $X$. Then, the operator norm of $\mathcal{G}$ is given by,
 \begin{equation}
\norm{\mathcal{G}}\ \ =\ \ \sup_{\norm{h}\, =\, 1} \norm{\mathcal{G} h}\ \ =\ \ 1.  
\end{equation}

\section{Two technical lemmas} 
\label{ttl}

\noindent 
Let $\mathfrak{C}_{n}$ denote the critical values of $T$ of order $n$, \textit{i.e.,} 
\[ \mathfrak{C}_{n}\ \ :=\ \ \bigcup_{k\, =\, 1}^{n} T^{k} \left( \mathfrak{C} \right). \] 
Consider a simply connected domain $U \subset X$ that satisfies $U_{\Sigma_{N}^{+}} = \Pi_{\Sigma_{N}^{+}} (U)$ being $\sigma$-invariant and $U_{\overline{\mathbb{C}}} = \Pi_{\overline{\mathbb{C}}} (U)$ remaining away from $\Pi_{\overline{\mathbb{C}}} \left( \mathfrak{C}_{n} \right)$, \textit{i.e.,} $\Pi_{\overline{\mathbb{C}}} \left( U \cap \mathfrak{C}_{n} \right) = \varnothing$. Hence, for any $m \le n$, there exists $(d_{1} + d_{2} + \cdots + d_{N})^{m}$ inverse branches of $T^{m}$, each analytic in $U$. 
\medskip 

\noindent 
\begin{lemma} 
\label{noofanalyticinverse}
Let $\gamma_{m} = \gamma_{m} (U)$ denote the total number of analytic inverse branches of $T^{m}$ in a simply connected domain $U \subset X$, and $\tau \le 2 \left( d_{1} + d_{2} + \cdots + d_{N} \right) - 2$ be the number of distinct critical values of $T$ in $U$. Then, 
\begin{equation} 
\label{inequilty1}
\left( d_{1} + d_{2} + \cdots + d_{N} \right)^{m} - \gamma_{m}\ \ \le\ \ \tau \sum_{t\, =\, 1}^{m\, -\, n} \left( d_{1} + d_{2} + \cdots + d_{N} \right)^{t}.
\end{equation}
\end{lemma} 
\medskip 

\noindent 
\begin{proof}(by induction) We start with the observation that the inequality in \eqref{inequilty1} is trivially true when $m \le n$. By inductive hypothesis, we assume it to be true for some $m > n$ and prove the same for $m + 1$. 
\medskip 

\noindent 
Since $U_{\overline{\mathbb{C}}}$ is simply connected, the Riemann-Hurwitz formula ensures that the $\gamma_{m}$ analytic inverse branches of $\Pi_{\overline{\mathbb{C}}} \left( T^{m} (U) \right)$ are mutually disjoint and simply connected. As there are $\tau$ distinct critical values of $T$ in $U$, at least $(\gamma_{m} - \tau)$ of these domains do not contain any critical values and each of these domains contain $(d_{1} + d_{2} + \cdots + d_{N})$ analytic inverse branches of the function $T$. Therefore we have at least $(d_{1} + d_{2} + \cdots + d_{N}) (\gamma_{m} - \tau)$ analytic inverse branches of $T^{m + 1}$ in $U$ (which is defined to be $\gamma_{m + 1}$). Thus,
\begin{eqnarray*}  
\gamma_{m + 1} &\ge & \left( d_{1} + d_{2} + \cdots + d_{N} \right) \left( \gamma_{m} - \tau \right) \\ 
& \ge & \left( d_{1} + d_{2} + \cdots + d_{N} \right)^{m + 1} - \tau \sum_{t\, =\, 1}^{m\, +\, 1\, -\, n} \left( d_{1} + d_{2} + \cdots + d_{N} \right)^{t}. 
\end{eqnarray*}
\end{proof}
\medskip 

\noindent 
\begin{lemma} 
\label{equicontinuous}
The family $\left\{ \mathcal{G}^{m} h : m \in \mathbb{Z}_{+} \right\}$ is equicontinuous for every function $h \in \mathcal{C} \left( \widetilde{\mathcal{J}}, \mathbb{R} \right)$. 
\end{lemma}  
\medskip 

\noindent  
\begin{proof} Let $h \ne 0$. Then $\forall\, \epsilon > 0,\ \exists\, \delta > 0$ such that 
\[ d_{X} \left((w_{1}, z_{1}), (w_{2}, z_{2})\right) < \delta\ \ \Longrightarrow\ \ d_{X} \left(h((w_{1}, z_{1})), h((w_{2}, z_{2}))\right) < \epsilon\ \ \forall\, (w_{i}, z_{i}) \in \widetilde{\mathcal{J}};\ i = 1, 2. \] 
In order to prove the lemma, we need to prove the equicontinuity of the family of sequence $\left\{ \mathcal{G}^{n} h \equiv h_{n} \right\}_{n\, \ge\, 1}$ in some $\epsilon$-neighbourhood of some generic point, say $(w_{0}, z_{0}) \in \widetilde{\mathcal{J}}$. We divide the arguments based on whether the point $(w_{0}, z_{0})$ belongs to the forward orbit of some critical point or not. 

\subsection*{Case 1: $(w_{0}, z_{0}) \notin \mathfrak{C}_{n}\ \forall n$} 
\label{case1}

\noindent 
Let $\epsilon > 0$. Choose $n$ large enough so that 
\[ \frac{4 \tau \norm{h}}{\left( d_{1} + d_{2} + \cdots + d_{N} \right)^{n}}\ \ <\ \ \epsilon. \] 
Choose $\alpha > 0$ in order that the neighbourhood 
\[ U_{\alpha}\ \ =\ \ \left\{ (w, z) \in X : d_{X} \left( (w, z), (w_{0}, z_{0}) \right) < \alpha \right\} \] 
satisfies $U_{\alpha} \cap \mathfrak{C}_{n} = \varnothing$. Then by lemma \eqref{noofanalyticinverse}, we have the estimate for the total number of analytic inverse branches of $T^{m}$ in $U_{\alpha}$. 
\begin{eqnarray*} 
\gamma_{m} & \ge & \left( d_{1} + d_{2} + \cdots + d_{N} \right)^{m} - \tau \sum_{t\, =\, 1}^{m\, -\, n} \left( d_{1} + d_{2} + \cdots + d_{N} \right)^{t} \\ 
& \ge & \left( d_{1} + d_{2} + \cdots + d_{N} \right)^{m} - 2 \tau \left( d_{1} + d_{2} + \cdots + d_{N} \right)^{m - n} \\
& \ge & \left( d_{1} + d_{2} + \cdots + d_{N} \right)^{m} \left( 1 - \frac{2 \tau}{\left( d_{1} + d_{2} + \cdots + d_{N} \right)^{n}} \right) \\ 
& > & \left( d_{1} + d_{2} + \cdots + d_{N} \right)^{m} \left( 1 - \frac{\epsilon}{2 \norm{h}} \right). 
\end{eqnarray*} 
This family of functions is normal, and hence equicontinuous. Therefore, there exists some $0 < \rho < \frac{\alpha}{2}$ such that for all $d_{X} \left( (w, z), (\omega, \zeta) \right) < \rho$, we have 
\[ d_{X} \left( T_{j}^{-m} \left((w, z)\right), T_{j}^{-m} \left((\omega, \zeta)\right) \right) < \delta,\ \forall\ (w, z)\ \text{and}\ (\omega, \zeta) \in U_{\frac{\alpha}{2}}\ \ \text{and}\ \ m \in \mathbb{Z}_{+}. \] 
Hence, taking $(w, z), (\omega, \zeta) \in \widetilde{\mathcal{J}}$ that are atmost $\rho$ apart from each other, we have, for $m$ sufficiently large 
\begin{eqnarray*} 
& & \left| h_{m} \left((w, z)\right) - h_{m} \left((\omega, \zeta)\right) \right| \ = \ \left| \mathcal{G}^{m} h \left((w, z)\right) - \mathcal{G}^{m} h \left((\omega, \zeta)\right) \right| \\ 
& = & \bigg| \frac{1}{\left( d_{1} + d_{2} + \cdots + d_{N} \right)^{m}} \sum_{v_{m}\, =\, 1}^{N} \cdots \sum_{v_{1}\, =\, 1}^{N} \sum_{j\, =\, 1}^{d_{v_{m}} \cdots d_{v_{1}}} h \left( (v_{m}\, \cdots\, v_{1}\, w, z_{((v_{m}\, \cdots\, v_{1}),\, j)}) \right) \\ 
& & -\ \frac{1}{\left( d_{1} + d_{2} + \cdots + d_{N} \right)^{m}} \sum_{v_{m}\, =\, 1}^{N} \cdots \sum_{v_{1}\, =\, 1}^{N} \sum_{j\, =\, 1}^{d_{v_{m}} \cdots d_{v_{1}}} h \left( (v_{m}\, \cdots\, v_{1}\, \omega, \zeta_{((v_{m}\, \cdots\, v_{1}),\, j)}) \right) \bigg|. 
\end{eqnarray*} 

\noindent 
Since there are only $\gamma_{m}$ analytic inverse branches of $T^{m}$, we have 
\begin{eqnarray*} 
& & \left| h_{m} \left((w, z)\right) - h_{m} \left((\omega, \zeta)\right) \right| \\ 
& \le & \frac{1}{(d_{1} + d_{2} + \cdots + d_{N})^{m}} \sum_{j\, =\, 1}^{\gamma_{m}} \bigg| \left( h \circ T_{j}^{-m} \right) \left((w, z)\right) - \left( h \circ T_{j}^{-m} \right) \left((\omega, \zeta)\right) \bigg| \\
& & +\ \frac{1}{(d_{1} + d_{2} + \cdots + d_{N})^{m}} 2 \left[ \left( d_{1} + d_{2} + \cdots + d_{N} \right)^{m} - \gamma_{m} \right] \norm{h} \\
& \le & \frac{1}{(d_{1} + d_{2} + \cdots + d_{N})^{m}} \left[ \gamma_{m} \epsilon + \frac{2 \left( d_{1} + d_{2} + \cdots + d_{N} \right)^{m} \norm{h} \epsilon}{2 \norm{h}} \right] \\ 
& \le & 2 \epsilon. 
\end{eqnarray*}
This proves the equicontinuity of $\left\{ h_{m} \right\}_{m \ge 1}$ in a neighbourhood of $(w_{0}, z_{0}) \notin \mathfrak{C}_{n}$ for any $n \in \mathbb{Z}_{+}$. 

\subsection*{Case 2: $(w_{0}, z_{0}) \in \mathfrak{C}_{n}$ for some $n \in \mathbb{Z}_{+}$} 
\subsubsection*{Subcase a: $(w_{0}, z_{0})$ is not a periodic point of any period for $T$}

\noindent 
Since $(w_{0}, z_{0})$ is not periodic for any period, there exists some $k \ge 1$ such that $(w_{0}, z_{0}) \notin \mathfrak{C}_{n} \setminus \mathfrak{C}_{k}$ for all $n > k$. Further, for any $l \ge 1,\ T^{l} \left((v_{n}\, \cdots\, v_{1}\, w_{0}), z_{0_{\left((v_{n}\, \cdots\, v_{1}),\, j\right)}}\right) \ne (w_{0}, z_{0})$ since $(w_{0}, z_{0})$ is not periodic. Moreover, for $n > k,\ \left((v_{n}\, \cdots\, v_{1}\, w_{0}), z_{0_{\left((v_{n}\, \cdots\, v_{1}),\, j\right)}}\right) \notin \mathfrak{C}_{q}$ for any $q$. For otherwise, there exists a critical point $(\omega_{c}, \zeta_{c}) \in \mathfrak{C}$ such that $ T^{q}(\omega_{c}, \zeta_{c}) = \left((v_{n}\, \cdots\, v_{1}\, w_{0}), z_{0_{\left((v_{n}\, \cdots\, v_{1}),\, j\right)}}\right)$ for some $q \in \mathbb{Z}_{+}$. Then $(\omega_{c}, \zeta_{c}) = T^{-\, (n + q)} (w_{0}, z_{0})$ resulting in $n + q \le k$ (from the first sentence), a contradiction to $k \le n$. 
\medskip 

\noindent 
Choose $\rho > 0$ small enough that the neighbourhoods, $U_{((v_{n}\, \cdots\, v_{1}),\, j)}$ around the $n$-th order pre-images of $(w_{0}, z_{0})$ in $\widetilde{\mathcal{J}}$ remain mutually disjoint, \textit{i.e.}, 
\begin{eqnarray*} 
U_{\left((v_{n}\, \cdots\, v_{1}),\, j\right)} & := & \left\{ (w, z) \in X : d_{X} \left( (w, z), \left((v_{n}\, \cdots\, v_{1}\, w_{0}), z_{0_{\left((v_{n}\, \cdots\, v_{1}),\, j\right)}}\right) \right) < \rho \right\}\ \cap\ \widetilde{\mathcal{J}}, \\ 
& & U_{\left((v_{n}\, \cdots\, v_{1}),\, j\right)} \cap U_{\left((\vartheta_{n}\, \cdots\, \vartheta_{1}),\, k\right)}\ \ \ne\ \ \varnothing\ \text{iff}\ v_{i} \equiv \vartheta_{i}\ \forall 1 \le i \le n\ \ \text{and}\ \ j \equiv k. 
\end{eqnarray*}

\noindent 
Since $\left((v_{n}\, \cdots\, v_{1}\, w_{0}), z_{0_{\left((v_{n}\, \cdots\, v_{1}),\, j\right)}}\right) \notin \mathfrak{C}_{q}$, the sequence $\left\{h_{m}\right\}$ is equicontinuous in each of these neighbourhoods, $U_{\left((v_{n}\, \cdots\, v_{1}),\, j\right)}$'s, by the preceding case. This means, for a fixed $\epsilon > 0$, there exists a $\delta > 0$ such that for all $m \in \mathbb{Z}_{+}$, we have 
\[ \left| h_{m} \left((w, z)\right) - h_{m} \left((\omega, \zeta)\right) \right| < \epsilon\ \ \text{if}\ \ d_{X}\left(( w, z), (\omega, \zeta)\right) < \delta\ \ \forall\ (w, z),\ (\omega, \zeta) \in U_{\left((v_{n}\, \cdots\, v_{1}),\, j\right)}. \] 

\noindent 
Choose $\eta > 0$ small enough such that whenever $d_{X} \left((w_{0}, z_{0}), (\omega_{0}, \zeta_{0})\right) < \eta$ with $(\omega_{0}, \zeta_{0}) \in \widetilde{\mathcal{J}}$, we have  
\[ d_{X} \left( \left((v_{n}\, \cdots\, v_{1}\, w_{0}), z_{0_{\left((v_{n}\, \cdots\, v_{1}),\, j\right)}}\right), \left((v_{n}\, \cdots\, v_{1}\, \omega_{0}), \zeta_{0_{\left((v_{n}\, \cdots\, v_{1}),\, j\right)}}\right) \right)\ \ <\ \ \rho. \]
Then, observe that 
\begin{eqnarray*} 
h_{n + m} (w, z) & = & \mathcal{G}^{n + m}(h)(w, z) = \mathcal{G}^{n} \left(\mathcal{G}^{m}(h)\right)(w, z) \\ 
& = & \frac{1}{(d_{1} + d_{2} + \cdots + d_{N})^{n}} \sum_{T^{n} (w_{n},\, z_{n})\, =\, (w,\, z)} \mathcal{G}^{m}(h)(w_{n}, z_{n}). 
\end{eqnarray*} 
Thus, for $d_{X} \left((w, z), (\omega, \zeta)\right) < \delta$ that satisfies 
\[ d_{X} \left((w, z), (w_{0}, z_{0})\right) < \eta\ \ \ \text{and}\ \ \ d_{X} \left((\omega, \zeta), (w_{0}, z_{0})\right) < \eta, \] 
we have 
\begin{eqnarray*} 
& & \left| h_{m + n}(w, z) - h_{m + n}(\omega, \zeta) \right| \\ 
& \le & \frac{1}{(d_{1} + d_{2} + \cdots + d_{N})^{n}} \\ 
& & \times\ \left[ \sum_{v_{n}\, =\, 1}^{N} \cdots \sum_{v_{1}\, =\, 1}^{N} \sum_{j\, =\, 1}^{d_{v_{n}} \cdots d_{v_{1}}} \left| h_{m} \left((v_{n}\, \cdots\, v_{1}\, w), z_{\left((v_{n}\, \cdots\, v_{1}),\, j\right)}\right) - h_{m} \left((v_{n}\, \cdots\, v_{1}\, \omega), \zeta_{\left((v_{n}\, \cdots\, v_{1}),\, j\right)}\right) \right| \right]. 
\end{eqnarray*} 
Suppose the points $(w, z)$ and $(\omega, \zeta)$ are atmost $\eta$ apart from each other, then it follows that for all $m > n,\ \left| h_{m}(w, z) - h_{m}(\omega, \zeta) \right| < \epsilon$. This is also true for $m \le n$ provided $(w, z)$ and $(\omega, \zeta)$ are chosen sufficiently close to each other. 

\subsubsection*{Subcase b: $T^{p} \left((w_{0}, z_{0})\right) = (w_{0}, z_{0})$ for some least $p \in \mathbb{Z}_{+}$} 

\noindent 
Since $ (w_{0}, z_{0})$ is a periodic point of period $p$, we have for all $n \in \mathbb{Z}_{+},\ \left( T^{pn} \right)' \left((w_{0}, z_{0})\right) \ne 0$, lest $(w_{0}, z_{0})$ violates our assumption (A1) by becoming a super-attracting cycle. 
\medskip 

\noindent 
For any $n \ge 1$, consider the pre-images of order $pn$ given by 
\[ \bigcup_{v_{pn}\, =\, 1}^{N} \cdots \bigcup_{v_{1}\, =\, 1}^{N} \bigcup_{j\, =\, 1}^{d_{v_{pn}} \cdots d_{v_{1}}} \left( \left((v_{pn}\, \cdots\, v_{1}\, w_0), z_{0_{((v_{pn}\, \cdots\, v_{1}),\, j)}}\right) \right). \] 
Note that among all these $(d_{1} + \cdots + d_{N})^{pn}$ pre-images, all except one point are different from $(w_{0}, z_{0})$. However, owing to periodicity, there exists one pre-image of order $pn$ that is equal to $(w_{0}, z_{0})$. Thus, this point  is necessarily of the form 
\[ v_{pn}\, \cdots\, v_{1}\, w_{0} = \underbrace{(w_{0})_{1}\, \cdots\, (w_{0})_{p}}_{n\ \text{times}} w_{0}\ \text{and the corresponding}\ z_{0_{((v_{pn}\, \cdots\, v_{1}),\, j)}} = z_{0},\ \text{for one}\ j. \]
Thus, any pre-image $\left((v_{pn}\, \cdots\, v_{1}\, w_{0}), z_{0_{((v_{pn}\, \cdots\, v_{1}),\, j)}} \right)$ other than $(w_{0}, z_{0})$ is not periodic. Suppose the point $\left((v_{pn}\, \cdots\, v_{1}\, w_{0}), z_{0_{((v_{pn}\, \cdots\, v_{1}),\, j)}}\right) \notin \mathfrak{C}_{q}$ for any $q$, then the arguments in Case 1 give the equicontinuity of the sequence $\left\{h_{m}\right\}$ in some neighbourhood of each such point. And whenever the pre-image $\left((v_{pn}\, \cdots\, v_{1}\, w_{0}), z_{0_{((v_{pn}\, \cdots\, v_{1}),\, j)}}\right) \in \mathfrak{C}_{q}$, the arguments in Subcase a of Case 2 yield the sequence $\left\{h_{m}\right\}$ to be equicontinuous in some neighbourhood of every such point. 
\medskip 

\noindent 
Now, for any $\epsilon > 0$, choose $n$ large enough such that $\frac{2 \norm{h}}{(d_{1} + \cdots + d_{N})^{pn}} < \epsilon$. Then there exists $\delta > 0$ and $\eta > 0$ so that for $(w, z), (\omega, \zeta) \in \widetilde{\mathcal{J}}$ satisfying 
\[ d_{X} \left((w, z), (\omega, \zeta)\right) < \delta,\ \ \ d_{X} \left((w, z), (w_{0}, z_{0})\right) < \eta,\ \ \ d_{X} \left((\omega, \zeta), (w_{0}, z_{0})\right) < \eta, \] 
we have for all $m \in \mathbb{Z}_{+}$ 
\begin{eqnarray*} 
& & \left| h_{pn\, +\, m} (w, z) - h_{pn\, +\, m} (\omega, \zeta) \right| \\ 
& \le & \frac{1}{(d_{1} + \cdots + d_{N})^{pn}} \\ 
& & \times\ \left[ \sum_{v_{pn}\, =\, 1}^{N} \cdots \sum_{v_{1}\, =\, 1}^{N} \sum_{j\, =\, 1}^{d_{v_{pn}} \cdots d_{v_{1}}} \left| h_{m} \left((v_{pn}\, \cdots\, v_{1}\, w), z_{\left((v_{pn}\, \cdots\, v_{1}),\, j\right)}\right) - h_{m} \left((v_{pn}\, \cdots\, v_{1}\, \omega), \zeta_{\left((v_{pn}\, \cdots\, v_{1}),\, j\right)}\right) \right| \right] \\  
& \le & \frac{1}{(d_{1} + \cdots + d_{N})^{pn}} \left[ \left\{ \left( d_{1} + \cdots + d_{N} \right)^{pn} - 1 \right\} \epsilon + 2 \norm{h} \right] \\ 
& < & 2 \epsilon. 
\end{eqnarray*} 
So, for $m > pn$, we have $\left| h_{m}(w, z) - h_{m}(\omega, \zeta) \right| < 2 \epsilon$. Choosing $\delta$ sufficiently small, one proves that the statement of lemma \eqref{equicontinuous} holds true for all $m$. 
\end{proof}

\section{Proof of theorem \eqref{edpr}} 
\label{pfpt2} 

\noindent 
Theorem \eqref{edpr} is concerned with the distribution of pre-images of a generic point $(w, z) \in \widetilde{\mathcal{J}}$. However, just as in the case of a single rational map, the properties of $\widetilde{\mathcal{J}}$ permit us to follow the arguments presented in \cite{ns:93}. However, one should observe that the arguments are not merely \textit{mutatis mutandis}, but a work towards accommodating the skew-product map. 
\medskip 

\begin{proof} We begin this proof with the following observation. For any continuous function $h : \widetilde{\mathcal{J}} \longrightarrow \mathbb{R}$, we have 
\begin{eqnarray*} 
\norm{\mathcal{G}^{m} h} & = & \sup_{(w,\, z)\, \in\, \widetilde{\mathcal{J}}} \left| \left( \mathcal{G}^{m} h \right) (w, z) \right| \\ 
& = & \sup_{(w,\, z)\, \in\, \widetilde{\mathcal{J}}} \left[ \frac{1}{\left( d_{1} + d_{2} + \cdots + d_{N} \right)^{m}} \right] \\ 
& & \times\ \left[ \sum_{v_{m}\, =\, 1}^{N} \cdots \sum_{v_{1}\, =\, 1}^{N} \sum_{j\, =\, 1}^{d_{v_{m}} \cdots d_{v_{1}}} \left| h \left( \left( \left(v_{m}\, \cdots\, v_{1}\, w \right), z_{\left((v_{m}\, \cdots\, v_{1}),\, j\right)} \right) \right) \right| \right] \\ 
& \le & \norm{h}. 
\end{eqnarray*} 
Thus, 
\[ \norm{\mathcal{G}^{m} h} = \norm{h}\ \ \text{iff}\ \ h \left( \left( \left(v_{m}\, \cdots\, v_{1}\, w \right), z_{\left((v_{m}\, \cdots\, v_{1}),\, j\right)} \right) \right) = \lambda \norm{h}. \] 
When $m$ is sufficiently large, the $\epsilon$-neighbourhood around each of these $m$-th order pre-images $\left( \left(v_{m}\, \cdots\, v_{1}\, w \right), z_{\left((v_{m}\, \cdots\, v_{1}),\, j\right)} \right)$ cover $\widetilde{\mathcal{J}}$. Hence, we can have only one of the following two cases; 
\begin{enumerate} 
\item[(a)] either $\norm{\mathcal{G}^{m} h} < \norm{h}$ for some $m$, and hence for every larger $m$; 
\item[(b)] or $\norm{\mathcal{G}^{m} h} = \norm{h}$ for all $m$, thereby making $h$ a constant.
\end{enumerate} 
\medskip

\noindent
Either way, this implies that $\left\{ \mathcal{G}^{m} h \right\}$ is uniformly bounded. Thus, by lemma \eqref{equicontinuous}, it is clear that the sequence, $\left\{ \mathcal{G}^{m} h \right\}$ is equicontinuous. An application of the Arzela-Ascoli theorem (as in \cite{gv:61}) then yields a subsequence, say $\left\{ \mathcal{G}^{m_{k}} h \right\}$ that uniformly converges to a limit $h_{\infty}$. Suppose $n_{k} = m_{k + 1} - m_{k}$. Then $m_{k} \to \infty$ implies $n_{k} \to \infty$. Consider the sequence $\left\{ \mathcal{G}^{n_{k}} h_{\infty} \right\}$ that uniformly converges to $H$. Then, 
\[ \norm{H - h_{\infty}}\ \ =\ \ \lim_{m_{k} \to \infty} \norm{\mathcal{G}^{m_{k + 1}} h_{\infty} - \mathcal{G}^{n_{k}} h}\ \ =\ \ \lim_{m_{k} \to \infty} \norm{\mathcal{G}^{n_{k}} \left( \mathcal{G}^{m_{k}} h_{\infty} - h \right)}\ \ =\ \ 0. \] 
Further, 
\[ \norm{H}\ \ =\ \ \lim_{k \to \infty} \norm{\mathcal{G}^{n_{k}} H}\ \ =\ \ \lim_{k \to \infty} \norm{\mathcal{G}^{n_{k} - m} \left( \mathcal{G}^{m} H \right)}\ \ \le\ \ \norm{\mathcal{G}^{m} H},\ \forall\, m \in \mathbb{Z}_{+}. \]
Therefore, the function $H = L_{h}$ is a constant that obeys 
\[ \norm{\mathcal{G}^{m} h - L_{h}}\ \ =\ \ \norm{\mathcal{G}^{m - m_{k}} \left( \mathcal{G}^{m_{k}} \left( h - L_{h} \right) \right)}\ \ \le\ \ \norm{\mathcal{G}^{m_{k}} \left( h - L_{h} \right)}. \] 
This implies that $\mathcal{G}^{m} h \to L_{h}$ as $m \to \infty$. Also 
\[ \norm{L_{h} - L_{h^{*}}}\ \ =\ \ \lim_{n \to \infty} \norm{\mathcal{G}^{n} \left( h - h^{*} \right)}\ \ \le\ \ \norm{h - h^{*}}. \]
Hence $L$ may be regarded as functional on the space of continuous functions on $\widetilde{\mathcal{J}}$. Since $L$ is bounded, linear and has norm $1$, Riesz representation theorem gives us a unique Borel probability measure $\mu^{*}$ such that, 
\[ L_{h} = \int_{\widetilde{\mathcal{J}}} h \left((\omega, \zeta)\right) d \mu^{*}(\omega, \zeta) \] 
and by construction, $\mu_{n}^{(w,\, z)} \rightharpoonup \mu^{*}$ independent of $(w, z)$.
\end{proof}

\section{Proof of theorem \eqref{lyuper}}
\label{pfper} 

\noindent 
In this section, we shall prove that the sequence of probability mass distributions $\left\{ \nu_{m} \right\}$ supported on the periodic points of order $m$ in $X$ converges to the measure $\mu^{*}$ obtained in the last section, in the weak*-topology, . The proof is based on the density of periodic points of $T$ in its Julia set $\widetilde{\mathcal{J}}$. Before we begin the proof of theorem \eqref{lyuper}, we make an observation about the set of critical points $\mathfrak{C}$ of $T$.  
\medskip 

\noindent 
Even though every map $R_{i}$ has only finitely many critical points for $1 \le i \le N$, one can observe that the skew-product map $T$ has uncountably many critical points. However for purposes of this proof, we will only be interested in the set of points $\Pi_{\overline{\mathbb{C}}} \left(\mathfrak{C}\right)$, that is finite. Here, we are interested in constructing a connected set $U \subset X$ of full measure $\mu^{*}$.
\medskip 

\noindent 
\begin{proof} 
Let $U \subset X$ be some open, connected set. Fix $\epsilon > 0$ and choose $n \in \mathbb{Z}_{+}$ such that $U$ remains away from $\mathfrak{C}_{n}$ and $T^{m}$ has at least $(d_{1} + \cdots + d_{N})^{m} (1 - \epsilon)$ number of analytic inverse functions in $U$. 
\medskip 

\noindent 
Now for any point $(\omega, \zeta) \in \mathfrak{C}_{n}$, let $r_{\zeta}$ denote the ray from $\zeta$ in $\overline{\mathbb{C}}$ to the point at $\infty$. The rays corresponding to different points in $\Pi_{\overline{\mathbb{C}}} \left(\mathfrak{C}_{n}\right)$ are mutually disjoint, except for their common terminal point. Since $\Pi_{\overline{\mathbb{C}}} \left(\mathfrak{C}_{n}\right)$ is atmost a finite set, we will only have atmost finitely many such rays in $\overline{\mathbb{C}}$. If $\mu^{*} \left( \left\{ \omega \right\} \times r_{\zeta} \right) > 0$ for some $(\omega, \zeta)$, then replace $r_{\zeta}$ by a new ray (also denoted with the same notation) obtained by rotating $r_{\zeta}$ through some small angle around $\zeta$ in $\overline{\mathbb{C}}$, in order that it avoids intersecting the other rays so that $\mu^{*} \left( \left\{ \omega \right\} \times r_{\zeta} \right) = 0$. 
\medskip 

\noindent 
Define $U$ such that $\Pi_{\overline{\mathbb{C}}} (U)$ is the entire Riemann sphere cut along these rays $r_{\zeta}$, \textit{i.e.}, 
\[ U\ \ :=\ \ \Sigma_{N}^{+} \times \left( \overline{\mathbb{C}} \setminus \bigcup\limits_{(\omega,\, \zeta)\, \in\, \mathfrak{C}_{n}} r_{\zeta} \right). \] 

\noindent 
Fix an open, non-empty set $E \subset X$ such that $\Pi_{\Sigma_{N}^{+}} (E)$ is a cylinder set of length $s$ that occurs at the beginning of the infinite-lettered word, say  $\left[ v_{1}\, v_{2}\, \cdots\, v_{s} \right],\ \Pi_{\overline{\mathbb{C}}} (E)$ is an open set in $\overline{\mathbb{C}}$ with $\mu^{*} (\partial E) = 0 $. Choose open sets $V$ and $W$ in $X$ such that $V$ is compactly contained in $U \cap E$ satisfying $\mu^{*}(V) > \mu^{*}(E) - \epsilon$ with the inclusions $\overline{V} \subset W \subset \overline{W} \subset U \cap E$ and $\Pi_{\Sigma_{N}^{+}} (V) = \Pi_{\Sigma_{N}^{+}} (W) = \Pi_{\Sigma_{N}^{+}} (E)$. 
\medskip 

\noindent 
Suppose the set $\Pi_{\overline{\mathbb{C}}} (V)$ is separated by a distance of $2l$ from $\Pi_{\overline{\mathbb{C}}} (\partial W)$. Denote by $T_{j}^{-m}$, the analytic inverse functions of $T^{m}$ in $U$, where $1 \le j \le \gamma_{m}$. Then, $\gamma_{m} \ge (d_{1} + \cdots + d_{N})^{m} (1 - \epsilon)$ is the number of analytic inverse functions of $T^{m}$ in $U$. This family of inverse functions is normal and thus, their limit function is constant on the Julia set $\widetilde{\mathcal{J}}$. Therefore we can choose $m > m_{0}$ such that 
\begin{equation} 
\label{diam}
{\rm diam} \left( \Pi_{\overline{\mathbb{C}}} \left(T_{j} ^{-m}\right) \left(\overline{W}\right) \right)\ <\ l, \quad \text{for} \quad  1 \le j \le \gamma_{m}. 
\end{equation}

\noindent 
Now consider a point $(\omega, \zeta) \in U \cap \widetilde{\mathcal{J}}$. Say that it has $n_{m} \le \gamma_{m}$ regular pre-images in $V$ given by $(w_{j}^{m}, z_{j}^{m}) \in T_{j}^{-m} (\omega, \zeta)$ for $1 \le j \le n_{m}$ and a maximum of $(d_{1} + \cdots + d_{N})^{m} - \gamma_{m} \le \epsilon (d_{1} + \cdots + d_{N})^{m}$ number of other pre-images. Thus we get,
\[ \mu_{m}^{(\omega,\, \zeta)} (V)\ \ <\ \ \frac{n_{m}}{(d_{1} + \cdots + d_{N})^{m}} + \epsilon. \] 
For $1 \le j \le n_{m}$, define the sets 
\[ U_{j}^{m}\ \ :=\ \ \left[ (w_{j}^{m})_{1}\, (w_{j}^{m})_{2}\, \cdots\, (w_{j}^{m})_{m} \right] \times \left\{ z \in \overline{\mathbb{C}} : d_{\overline{\mathbb{C}}} (z, z_{j}^{m}) < 2l \right\}. \]

\noindent 
As $(w_{j}^{m}, z_{j}^{m}) \in V$, we have by equation \eqref{diam}  
\[ T_{j}^{-m} \left(\overline{W}\right)\ \ \subset\ \ U_{j}^{m}\ \ \subset\ \ W,\ \ \ \ \text{for}\ \ m > m_{0}. \] 
In particular, $T_{j}^{-m} : \overline{U_{j}^{m}} \longrightarrow \overline{U_{j}^{m}}$ is a compact mapping. Hence, for every $1 \le j \le n_{m}$, the inverse function $T_{j}^{-m}$ has a fixed point say $(\omega_{j}^{m}, \zeta_{j}^{m})$ in $U_{j}^{m}$. Thus $(\omega_{j}^{m}, \zeta_{j}^{m})$ are fixed points of $T^{m}$ in $W$ and hence $T^{m}$ has at least $n_{m}$ distinct fixed points in $W$, since the image domains $T_{j}^{-m} (W) \subset W$ are disjoint. Hence, for sufficiently large $m$
\begin{eqnarray*} 
\nu_{m} (E) & \ge & \nu_{m} (W) \\ 
& \ge & \frac{n_{m}}{(d_{1} + \cdots + d_{N})^{m} + N} \\
& \ge & \mu_{m}^{(\omega,\, \zeta)} (U) \frac{(d_{1} + \cdots + d_{N})^{m}}{(d_{1} + \cdots + d_{N})^{m} + N} - \epsilon. 
\end{eqnarray*} 
Taking limit infimum as $m$ tends to infinity and bearing in mind that $\epsilon$ is an arbitrarily small positive quantity, we get, 
\[ \nu (E)\ \ \ge\ \ \liminf_{m \to \infty} \nu_{m} (E)\ \ \ge\ \ \mu^{*} (E). \] 
This is true for all Borel sets and since both of them are probability measures that equidistributes points, we have $\mu^{*}$ and $\nu$ agree on all Borel sets. 
\end{proof}

\section{Proofs of theorems \eqref{lard-preimage} and \eqref{periodicpts}} 
\label{pfld} 

\noindent 
We begin this section with a remark that the statements of the theorems \eqref{lard-preimage} and \eqref{periodicpts} are merely stating certain large deviation results, as can be compared with \cite{lsy:90, ps:96, ps:07}. 
\medskip 

\noindent 
Throughout this section, we take $(w, z) \in \widetilde{\mathcal{J}},\ f : \widetilde{\mathcal{J}} \longrightarrow \mathbb{R}$ to be a H\"older continuous function that satisfies $\mathfrak{P} (f) > \sup f + \log N$ and $g \in \mathcal{C} \left( \widetilde{\mathcal{J}}, \mathbb{R} \right)$. Let $\mathcal{M}$ denote the set of all probability measures supported on $\widetilde{\mathcal{J}}$ while we recall that $\mathcal{M}_{T} \subset \mathcal{M}$ denotes the set of all $T$-invariant probability measures supported on $\widetilde{\mathcal{J}}$.  
\medskip 

\noindent 
\begin{lemma}\cite{ps:96}
The map $m \longmapsto h_{m} (T)$ is upper semicontinuous in weak* topology on $\mathcal{M}_{T}$. Further, for $m \in \mathcal{M}_{T}$, we have 
\[ h_{m} (T)\ \ =\ \ \inf \left\{ \mathfrak{P} (g) - \int g d m\ :\ g \in \mathcal{C} \left( \widetilde{\mathcal{J}}, \mathbb{R} \right) \right\}. \] 
\end{lemma}
\medskip 

\noindent 
In fact, this lemma is true for any continuous mapping $T$ defined on a compact metric space that has finite topological entropy. In order to prove theorem \eqref{lard-preimage}, we now define a function $Q$ and its Legendre transform $\mathbb{L}$ by 
\begin{eqnarray*} 
Q : \mathcal{C} \left( \widetilde{\mathcal{J}}, \mathbb{R} \right) \longrightarrow \mathbb{R} & \text{given by} & Q (g) := \mathfrak{P} \left( f + g \right) - \mathfrak{P} \left( f \right), \\ 
\mathbb{L} : \mathcal{M}_{T} \longrightarrow \mathbb{R} & \text{given by} & \mathbb{L} \left( m \right) := \sup_{g\, \in\, \mathcal{C} \left( \widetilde{\mathcal{J}} \right)} \left\{ \int g d m - Q (g) \right\}. 
\end{eqnarray*} 

\noindent 
Further, for any weak*-closed subset $\mathcal{K} \subset \mathcal{M}_{T}$, we define 
\begin{equation} 
\label{rhok} 
\rho \left(\mathcal{K}\right) := \inf_{m\, \in\, \mathcal{K}} \mathbb{L} (m). 
\end{equation}  
Then, we have the following two lemmas, the first corresponding to the sequence of normalised proportion of orbital measure and the second corresponding to the sequence of normalised proportion of periodic orbital measure. 
\medskip 

\noindent 
\begin{lemma} 
\label{lem}
Let $(w, z) \in \widetilde{ \mathcal{J} }$ be a generic point whose $n$-th order pre-images are collected in  
\[ \bigcup_{v_{n} = 1}^{N} \cdots \bigcup_{v_{1} = 1}^{N} \bigcup_{j = 1}^{d_{v_{n}} \cdots d_{v_{1}}} \left( \left( (v_{n}\, \cdots\, v_{1}\, w), z_{((v_{n}\, \cdots\, v_{1}),\, j)} \right) \right). \] 
For any choice of $n$-lettered word $(v_{n}\, v_{n - 1}\, \cdots\, v_{1})$ and an appropriate choice of $z_{0}$ \textit{i.e.}, 
\[ \left( R_{v_{1}} \circ R_{v_{2}} \circ \cdots \circ R_{v_{n}} \right) \left(z_{0}\right) = z, \] 
recall the definition of a sequence of the normalised proportion of orbital measure given by 
\begin{eqnarray*} 
\mu_{\left( \left( (v_{n}\, \cdots\, v_{1}\, w),\, z_{0} \right) \right)} & = & \frac{1}{n} \bigg[ \delta_{\left( \left( (v_{n}\, \cdots\, v_{1}\, w),\, z_{0} \right) \right)} + \delta_{T \left( \left( \left( (v_{n}\, \cdots\, v_{1}\, w),\, z_{0} \right) \right) \right)} \\ 
& & \ \ \ \ \ \ \ \ \ \ \ \ \ \ \ \ \ \ \ \ \ \ \ \ + \cdots + \delta_{T^{n - 1} \left( \left( \left( (v_{n}\, \cdots\, v_{1}\, w),\, z_{0} \right) \right) \right)} \bigg] \\ 
& = & \frac{1}{n} \bigg[ \delta_{\left( \left( (v_{n}\, \cdots\, v_{1}\, w),\, z_{0} \right) \right)} + \delta_{\left( \left( (v_{n - 1}\, \cdots\, v_{1}\, w),\, R_{v_{n}} z_{0} \right) \right)} \\ 
& & \ \ \ \ \ \ \ \ \ \ \ \ \ \ \ \ \ \ \ \ \ \ \ \ + \cdots + \delta_{\left( \left( \left( (v_{1}\, w),\, R_{v_{2}} \circ \cdots \circ R_{v_{n}} \right) \left(z_{0}\right) \right) \right)} \bigg]  
\end{eqnarray*} 
Then, for any weak*-closed subset $\mathcal{K} \subset \mathcal{M}_{T}$, we have 
\[ \limsup_{n \to \infty} \frac{1}{n} \log \left[ \frac{\sum\limits_{\mu_{\left( \left( (v_{n}\, \cdots\, v_{1}\, w),\, z_{0} \right) \right)}\, \in\, \mathcal{K}} e^{f^{n} \left( (v_{n}\, \cdots\, v_{1}\, w),\, z_{0} \right)}}{\sum e^{f^{n} \left( (v_{n}\, \cdots\, v_{1}\, w),\, z_{0} \right)}} \right]\ \ \le\ \ - \rho \left(\mathcal{K}\right). \]
\end{lemma}
\medskip 

\noindent 
\begin{lemma} 
\label{lm} 
Let $(\omega, \zeta) \in \widetilde{\mathcal{J}} \cap {\rm Fix}_{m} (T)$, \textit{i.e.}, $(\omega, \zeta) = \left( \left( \vartheta_{1}\, \vartheta_{2}\, \cdots\, \vartheta_{m}\, \vartheta_{1}\, \cdots\, \vartheta_{m}\, \vartheta_{1}\, \cdots\, \vartheta_{m}\, \cdots \right), \zeta \right)$, where $\left( R_{\vartheta_{m}} \circ \cdots \circ R_{\vartheta_{1}} \right) \left(\zeta\right) = \zeta$. Recall the definition of a sequence of the normalised proportion of the periodic orbital measure given by 
\begin{eqnarray*} 
\nu_{(\omega,\, \zeta),\, m} & = & \frac{1}{m} \left[ \delta_{(\omega,\, \zeta)} + \delta_{T \left((\omega,\, \zeta)\right)} + \cdots + \delta_{T^{m - 1} \left((\omega,\, \zeta)\right)} \right] \\ 
& = & \frac{1}{m} \bigg[ \delta_{\left( \left( \vartheta_{1}\, \vartheta_{2}\, \cdots\, \vartheta_{m}\, \vartheta_{1}\, \cdots\, \vartheta_{m}\, \cdots \right),\, \zeta \right)} + \delta_{\left( \left( \vartheta_{2}\, \cdots\, \vartheta_{m}\, \vartheta_{1}\, \cdots\, \vartheta_{m}\, \cdots \right),\, R_{\vartheta_{1}} \zeta \right)} \\ 
& & \ \ \ \ \ \ \ \ \ \ \ \ \ \ \ \ \ \ \ \ \ \ \ \ \ \ \ \ \ \ \ \ \ \ + \cdots + \delta_{\left( \left( \vartheta_{m}\,\vartheta_{1}\, \cdots\, \vartheta_{m}\, \cdots \right),\, \left( R_{\vartheta_{m - 1}} \circ \cdots \circ R_{\vartheta_{1}} \right) \zeta \right)} \bigg]
\end{eqnarray*}  
Then, for any weak*-closed subset $\mathcal{K} \subset \mathcal{M}_{T}$, we have 
\[ \limsup_{m \to \infty} \frac{1}{m} \log \left[ \frac{\sum\limits_{\substack{(\omega,\, \zeta)\, \in\, \widetilde{\mathcal{J}} \cap {\rm Fix}_{m} (T) \\ \nu_{(\omega,\, \zeta),\, m}\, \in\, \mathcal{K}}} e^{f^{m} (\omega,\, \zeta)}}{\sum\limits_{(\omega,\, \zeta)\, \in\, \widetilde{\mathcal{J}} \cap {\rm Fix}_{m} (T)} e^{f^{m} \left(\omega,\, \zeta\right)}} \right]\ \ \le\ \ - \rho \left(\mathcal{K}\right). \]
\end{lemma}
\medskip 

\noindent 
Before we prove the lemmas, we remark that as a fallout of theorems \eqref{edpr} and \eqref{lyuper} or otherwise from \cite{hs:00, hs:01, hs:06, su:09}, the set of all pre-images of any generic point in $\widetilde{\mathcal{J}}$ and the set of all periodic points of all periods form a dense subset of $\widetilde{\mathcal{J}}$. A proof of the above lemmas \eqref{lem} and \eqref{lm} essentially exploits the alternate definition of pressure, as written after theorem \eqref{SumiUrbanski} and the statement of theorem \eqref{ps:96}. Since a complete proof of both the lemmas run along similar lines, we shall contend ourselves by proving lemma \eqref{lem} only, in this paper. 
\medskip 

\noindent 
\begin{proof}(of lemma \eqref{lem}) 
Recall from equation \eqref{rhok} that $\rho \left(\mathcal{K}\right) = \inf_{m\, \in\, \mathcal{K}} \mathbb{L} (m)$ for any weak*-closed subset $\mathcal{K} \subset \mathcal{M}_{T}$. Hence, by definition, it is true that for every $\epsilon > 0$ and for every $m \in \mathcal{K}$, 
\[ \rho \left(\mathcal{K}\right) - \epsilon\ \ <\ \ \mathbb{L} (m)\ \ =\ \ \sup_{g\, \in\, \mathcal{C} \left( \widetilde{\mathcal{J}} \right)} \left\{ \int g d m - Q(g) \right\}. \] 
In other words, for every $m \in \mathcal{K}$ there exists a $g \in \mathcal{C} \left( \widetilde{\mathcal{J}}, \mathbb{R} \right)$ such that 
\[ \rho \left(\mathcal{K}\right) - \epsilon\ \ <\ \ \int g d m - Q(g). \] 
Thus, $\mathcal{K} \subset \mathcal{M}_{T}$ is covered by 
\[ \mathcal{K}\ \ \subset\ \ \bigcup_{g\, \in\, \mathcal{C} \left( \widetilde{\mathcal{J}}, \mathbb{R} \right)} \left\{ m \in \mathcal{M}_{T} : \int g d m - Q(g) > \rho \left(\mathcal{K}\right) - \epsilon \right\}. \] 
Since $T$ is a continuous map defined on the compact metric space $\widetilde{\mathcal{J}}$, the space of $T$-invariant probability measures $\mathcal{M}_{T}$ is weak*-compact. Thus $\mathcal{K}$, a weak*-closed subset of $\mathcal{M}_{T}$ is also compact. Hence, a finite subcover is sufficient for $\mathcal{K}$. Therefore, 
\[ \mathcal{K}\ \ \subset\ \ \bigcup_{i\, =\, 1}^{s} \left\{ m \in \mathcal{M}_{T} : \int g_{i} d m - Q(g_{i}) > \rho \left(\mathcal{K}\right) - \epsilon \right\}. \]

\noindent 
Let $(w, z) \in \widetilde{\mathcal{J}}$ be a generic point. Now for any choice of $n$-lettered word $(v_{n}\, v_{n - 1}\, \cdots\, v_{1})$ and an appropriate choice of $z_{0}$ satisfying  
\[ \left( R_{v_{1}} \circ R_{v_{2}} \circ \cdots \circ R_{v_{n}} \right) \left(z_{0}\right) = z, \] 
note that, 
\begin{eqnarray*} 
& & \left\{ \left( \left( (v_{n}\, \cdots\, v_{1}\, w), z_{0} \right) \right) \in \widetilde{\mathcal{J}} : \mu_{\left( \left( (v_{n}\, \cdots\, v_{1}\, w),\, z_{0} \right) \right)} \in \mathcal{K} \right\} \\ 
& \subset & \bigcup_{i\, =\, 1}^{s} \left\{ \left( \left( (v_{n}\, \cdots\, v_{1}\, w), z_{0} \right) \right) \in \widetilde{\mathcal{J}} : \int g_{i} d \mu_{\left( \left( (v_{n}\, \cdots\, v_{1}\, w),\, z_{0} \right) \right)} - Q(g_{i}) > \rho \left(\mathcal{K}\right) - \epsilon \right\} \\ 
& = & \bigcup_{i\, =\, 1}^{s} \left\{ \left( \left( (v_{n}\, \cdots\, v_{1}\, w), z_{0} \right) \right) \in \widetilde{\mathcal{J}} : \frac{1}{n} \left[ g_{i}^{n} \left( \left( (v_{n}\, \cdots\, v_{1}\, w), z_{0} \right) \right) \right] > Q(g_{i}) + \rho \left(\mathcal{K}\right) - \epsilon \right\}. 
\end{eqnarray*} 

\noindent 
Thus, we obtain the following inequality, 
\begin{eqnarray*} 
& & \sum\limits_{\substack{ \left( (v_{n}\, \cdots\, v_{1}\, w),\, z_{0} \right)\, \in\, \widetilde{\mathcal{J}} \\  \mu_{\left( \left( (v_{n}\, \cdots\, v_{1}\, w),\, z_{0} \right) \right)}\, \in\, \mathcal{K}}} e^{f^{n} \left( (v_{n}\, \cdots\, v_{1}\, w),\, z_{0} \right)} \\ 
& \le & \sum\limits_{i\, =\, 1}^{s} \sum\limits_{\substack{ \left( (v_{n}\, \cdots\, v_{1}\, w),\, z_{0} \right)\, \in\, \widetilde{\mathcal{J}} \\ \left( \left( (v_{n}\, \cdots\, v_{1}\, w),\, z_{0} \right) \right)\, -\, Q(g)\, >\, n \left[ \rho \left(\mathcal{K}\right)\, -\, \epsilon \right]}} e^{- g_{i}^{n} \left( (v_{n}\, \cdots\, v_{1}\, w),\, z_{0} \right)}\ e^{f^{n} \left( (v_{n}\, \cdots\, v_{1}\, w),\, z_{0} \right) + g_{i}^{n} \left( (v_{n}\, \cdots\, v_{1}\, w),\, z_{0} \right)} \\
& \le & \sum\limits_{i\, =\, 1}^{s} \sum\limits_{\left( (v_{n}\, \cdots\, v_{1}\, w),\, z_{0} \right)\, \in\, \widetilde{\mathcal{J}}} e^{-n \left[ Q \left(g_{i}\right)\, +\, \left( \rho\left(\mathcal{K}\right)\, -\, \epsilon \right) \right]}\ e^{f^{n} \left( (v_{n}\, \cdots\, v_{1}\, w),\, z_{0} \right)\, +\, g_{i}^{n} \left( (v_{n}\, \cdots\, v_{1}\, w),\, z_{0} \right)} \\ 
& = & \sum\limits_{i\, =\, 1}^{s} e^{-n \left[ Q \left(g_{i}\right)\, +\, \left(\rho\left(\mathcal{K}\right)\, -\, \epsilon \right) \right]}\ \sum\limits_{\left( (v_{n}\, \cdots\, v_{1}\, w),\, z_{0} \right)\, \in\, \widetilde{\mathcal{J}}} e^{f^{n} \left( (v_{n}\, \cdots\, v_{1}\, w),\, z_{0} \right)}\ e^{g_{i}^{n} \left( (v_{n}\, \cdots\, v_{1}\, w),\, z_{0} \right)} \\ 
& \le & s\ \sup\limits_{1\, \le\, i\, \le\, s} \left\{ e^{-n \left[ Q\left(g_{i}\right)\, +\, \left( \rho\left(\mathcal{K}\right)\, -\,\epsilon \right) \right]}\ \sum\limits_{\left( (v_{n}\, \cdots\, v_{1}\, w),\, z_{0} \right)\, \in\, \widetilde{\mathcal{J}}} e^{f^{n} \left( (v_{n}\, \cdots\, v_{1}\, w),\, z_{0} \right)}\ e^{g_{i}^{n} \left( (v_{n}\, \cdots\, v_{1}\, w),\, z_{0} \right)} \right\}.
\end{eqnarray*}

\noindent 
Finally, we consider,
\begin{eqnarray*}
& & \limsup_{n \to \infty} \frac{1}{n} \log \left[ \frac{\sum\limits_{\substack{\left( (v_{n}\, \cdots\, v_{1}\, w),\, z_{0} \right)\, \in\, \widetilde{\mathcal{J}} \\ \mu_{\left( \left( (v_{n}\, \cdots\, v_{1}\, w),\, z_{0} \right) \right)}\, \in\, \mathcal{K}}} e^{f^{n} \left( (v_{n}\, \cdots\, v_{1}\, w),\, z_{0} \right)}}{\sum\limits_{\left( (v_{n}\, \cdots\, v_{1}\, w),\, z_{0} \right)\, \in\, \widetilde{\mathcal{J}}} e^{f^{n} \left( (v_{n}\, \cdots\, v_{1}\, w),\, z_{0} \right)}} \right]\\
& \le & \sup\limits_{1\, \le\, i\, \le\, s} \bigg[ - Q\left(g_{i}\right)\, -\, \rho\left(\mathcal{K}\right)\, +\, \epsilon \\ 
& & \hspace{+2cm} +\ \limsup_{n \to \infty} \frac{1}{n} \log \sum\limits_{\left( (v_{n}\, \cdots\, v_{1}\, w),\, z_{0} \right)\, \in\, \widetilde{\mathcal{J}}} e^{f^{n} \left( (v_{n}\, \cdots\, v_{1}\, w),\, z_{0} \right)}\ e^{g_{i}^{n} \left( (v_{n}\, \cdots\, v_{1}\, w),\, z_{0} \right)} \\ 
& & \hspace{+4cm} -\ \liminf_{n \to \infty} \frac{1}{n} \log \sum\limits_{\left( (v_{n}\, \cdots\, v_{1}\, w),\, z_{0} \right)\, \in\, \widetilde{\mathcal{J}}} e^{f^{n} \left( (v_{n}\, \cdots\, v_{1}\, w),\, z_{0} \right)} \bigg] \\
\end{eqnarray*} 
\begin{eqnarray*}
& \le & \sup\limits_{1\, \le\, i\, \le\, s} \left[ - Q\left(g_{i}\right) - \rho\left(\mathcal{K}\right) + \epsilon + \mathfrak{P} \left( f + g_{i} \right) - \mathfrak{P} \left(f\right) \right] \\ 
& = & - \rho\left(\mathcal{K}\right) + \epsilon.
\end{eqnarray*}
Since $\epsilon$ is an arbitrarily small positive quantity, the proof of this lemma is complete.
\end{proof}
\medskip 

\noindent 
As one may expect, we will use lemma \eqref{lem} to prove theorem \eqref{lard-preimage} and lemma \eqref{lm} to prove theorem \eqref{periodicpts}. We finally state a lemma (with a short proof for the readers' convenience) due to Pollicott and Sharp, \cite{ps:96} that clinches the proof of the theorem \eqref{lard-preimage}. It ensures that if $\mu_{f} \notin \mathcal{K}$, then $\rho\left(\mathcal{K}\right) > 0$. 
\medskip 

\noindent 
\begin{lemma}\cite{ps:96}
\label{lard_lemma_last} 
If $m \ne \mu_{f}$, then $\mathbb{L}(m) > 0$. The map $m \longmapsto \mathbb{L} (m)$ is lower semicontinuous on $\mathcal{M}_{T}$. 
\end{lemma}
\medskip 

\noindent 
\begin{proof}
Consider 
\begin{eqnarray*} 
\mathbb{L} (m) & = & \sup_{g\, \in\, \mathcal{C} \left( \widetilde{\mathcal{J}} \right)} \left\{ \int g d m - Q (g) \right\} \\ 
& = & \sup_{g\, \in\, \mathcal{C} \left( \widetilde{\mathcal{J}} \right)} \left\{ \int g d m - \mathfrak{P} (f + g) + \mathfrak{P} (f) \right\} \\ 
& = & \sup_{g\, \in\, \mathcal{C} \left( \widetilde{\mathcal{J}} \right)} \left\{ \int \left( g - f \right) d m - \mathfrak{P} (g) + \mathfrak{P} (f) \right\} \\ 
& = & \sup_{g\, \in\, \mathcal{C} \left( \widetilde{\mathcal{J}} \right)} \left\{ \int g d m - \mathfrak{P} (g) \right\} + \mathfrak{P} (f) - \int f d m \\ 
& = & - \inf_{g\, \in\, \mathcal{C} \left( \widetilde{\mathcal{J}} \right)} \left\{ \mathfrak{P} (g) - \int g d m \right\} + \mathfrak{P} (f) - \int f d m \\ 
& = & - h(m) + \mathfrak{P} (f) - \int f d m \\ 
& > & 0.
\end{eqnarray*} 
\end{proof} 
\medskip 

\noindent 
\begin{proof}(of theorem \eqref{lard-preimage}) Let $\mu_{f} \notin \mathcal{K}$. Then this implies $\rho\left(\mathcal{K}\right) > 0$. For any weak*-open neighbourhood $\mathcal{U} \subset \mathcal{M}_{T}$ of $\mu_{f}$, define $\mathcal{K} = \mathcal{M}_{T} \setminus \mathcal{U}$. Then, $\mathcal{K}$ is weak*-closed, and hence weak*-compact. The proof of the theorem is then complete by using lemma \eqref{lem}. 
\end{proof} 
\medskip 

\noindent 
The proof of theorem \eqref{periodicpts} follows exactly along the same lines, except for using lemma \eqref{lm}, in place of lemma \eqref{lem}. 
\bigskip \bigskip 

\bibliographystyle{plainnat}
\bibliography{mybib} 

\bigskip 
\bigskip 

\noindent 
\textsc{Shrihari Sridharan} \\ 
Indian Institute of Science Education and Research Thiruvananthapuram (IISER-TVM), \\ 
Maruthamala P.O., Vithura, Thiruvananthapuram, INDIA. PIN 695 551. \\
{\tt shrihari@iisertvm.ac.in}  
\bigskip \\ 

\noindent 
\textsc{Sharvari Neetin Tikekar} \\ 
Indian Institute of Science Education and Research Thiruvananthapuram (IISER-TVM), \\ 
Maruthamala P.O., Vithura, Thiruvananthapuram, INDIA. PIN 695 551. \\
{\tt sharvai.tikekar14@iisertvm.ac.in} 
\bigskip \\ 

\noindent 
\textsc{Atma Ram Tiwari} \\ 
Indian Institute of Science Education and Research Thiruvananthapuram (IISER-TVM), \\ 
Maruthamala P.O., Vithura, Thiruvananthapuram, INDIA. PIN 695 551. \\
{\tt artiwari15@iisertvm.ac.in} 

\end{document}